\documentclass{bica}
\usepackage{color,graphicx}
\def\newpic#1{%
\def\emline##1##2##3##4##5##6{%
\put(##1,##2){\special{em:point #1##3}}%
\put(##4,##5){\special{em:point #1##6}}%
\special{em:line #1##3,#1##6}}}
\newpic{}
\def\emline#1#2#3#4#5#6{%
\put(#1,#2){\special{em:moveto}}%
\put(#4,#5){\special{em:lineto}}}
\def\newpic#1{}

\newcommand{\N}{\hbox{\bf N}}
\newcommand{\Z}{\hbox{\bf Z}}
\newcommand{\x}{\overline}
\newauthor{%
Italo J. Dejter}{%
Italo J. Dejter}{%
University of Puerto Rico\\
Rio Piedras, PR 00936-8377}[%
italo.dejter@gmail.com]

\title{Rainbow tetrahedra in Cayley graphs}


\classnbr{05C15; 05C75; 05C62}

\begin{document}
\begin{abstract}
\noindent Let $\Gamma_n$ be the complete undirected Cayley graph of the odd
cyclic group $\Z_n$. Connected graphs whose vertices
are rainbow tetrahedra in $\Gamma_n$ are studied, with any two such vertices adjacent if and only if they share (as tetrahedra) precisely two distinct triangles.
This yields graphs $G$ of largest degree 6, asymptotic diameter $|V(G)|^{1/3}$ and almost all vertices with degree: {\bf(a)} 6 in $G$; {\bf(b)} 4 in exactly six connected subgraphs of the $(3,6,3,6)$-semi-regular tessellation; and {\bf(c)} 3 in exactly four connected subgraphs of the $\{6,3\}$-regular hexagonal tessellation. These vertices have as closed neighborhoods the union (in a fixed way) of closed neighborhoods in the ten respective resulting tessellations. 
\end{abstract}

\section{Introduction}\label{s1}

Cayley graphs are very important because they have many useful applications (cf. \cite{Kelarev1})
and are related to automata theory (cf. \cite{Kelarev2,Kelarev3}). In the present work, we deal with Cayley graphs of  a finite abelian group $G$ with its identity denoted 0. Let $S$ be a subset of $G$
such that $0\notin S$ and $S=-S$ (that is: $s\in S$ if and only if $-s\in S$).
The {\it Cayley graph} $\Gamma(G,S)$ on $G$ with {\it connection set} $S$ is a graph that has as its vertices  the elements of $G$ and is such that it has an edge $e$ joining vertices $g$ and $h$ if and only if
$h=g+s$, for some $s\in S$. In this case, we say that the edge $e$ has {\it color} $s$.
A concept of ``rainbow'' has been used in various fashions in a graph theory context, in
\cite{Aha,Alon,Bar,Fri,HMP,Jan,Kos,LeS,MS,OYY,Per,SR,Wang,Woldar} and related papers. Ours is in relation to edge colors in Cayley graphs of finite cyclic groups.
Below, the complete graph $K_n=K_{2k+1}$ will be viewed as the
Cayley graph $\Gamma_n = \Gamma(\Z_n,[k])$ of the cyclic group $\Z_n$ of integers mod $n$ with connecting set $[k]=$ $\{1,2,...,$ $k\}$.
Relations among {\it rainbow} triangles and tetrahedra  in $\Gamma_n$ ({\it rainbow} meaning here edges with pairwise different colors) will be shown to yield a family ${\mathcal G}_1$ of
connected graphs $G=G_{n,4}$ of largest degree $\Delta(G)=6$, asymptotic
diameter $|V(G)|^{1/3}$ and such that
almost all its vertices $v$ have degree: {\bf(a)} 6 in $G$; {\bf(b)} 4 in exactly six connected subgraphs of the $(3,6,3,6)$-semi-regular tessellation (\cite{Fej}, page 43); and {\bf(c)} 3 in exactly four connected subgraphs of the $\{6,3\}$-regular hexagonal tessellation (\cite{Fej}, page 43).
We refer to each of these ten subgraphs of $G$ as a $\mathcal D$- or as an $\mathcal H$-{\it modeled subgraph} of $G$ if it is as in (b) or as in (c) above, respectively. On the other hand, based on rainbow triangles a family ${\mathcal G}_0$ of connected graphs $G=G_{n,3}$ of largest degree $\Delta(G)=3$ and asymptotic
diameter $|V(G)|^{1/2}$  was introduced in \cite{DHS}. See Section~\ref{s3} below for a short survey of \cite{DHS} and for further developments ahead in this paper.

The mentioned asymptotic properties of the families ${\mathcal G}_0$ and ${\mathcal G}_1$ confirm the following conjecture, further discussed in \cite{appendix}.

\begin{con}\label{conj3}
The asymptotic diameter of a family of graphs $G$ with a common $\Delta(G)$ is a given $($radical, logarithmic, $\ldots)$ function of the vertex number of $G$.
\end{con}

\section{Main results}\label{s2}

The present paper is devoted to the following results, containing the claimed properties of ${\mathcal G}_1$. (For related properties, see \cite{De4} and its references).
The {\it tessellated neighborhood} of a vertex $v$ in a $\mathcal D$- or $\mathcal H$-modeled subgraph $G$
is formed by $v$ and its incident edges and faces as well as by the other edges adjacent to those faces and the endvertices of these edges.

\begin{figure}[htp]
\hspace*{1mm}
\includegraphics[scale=0.29]{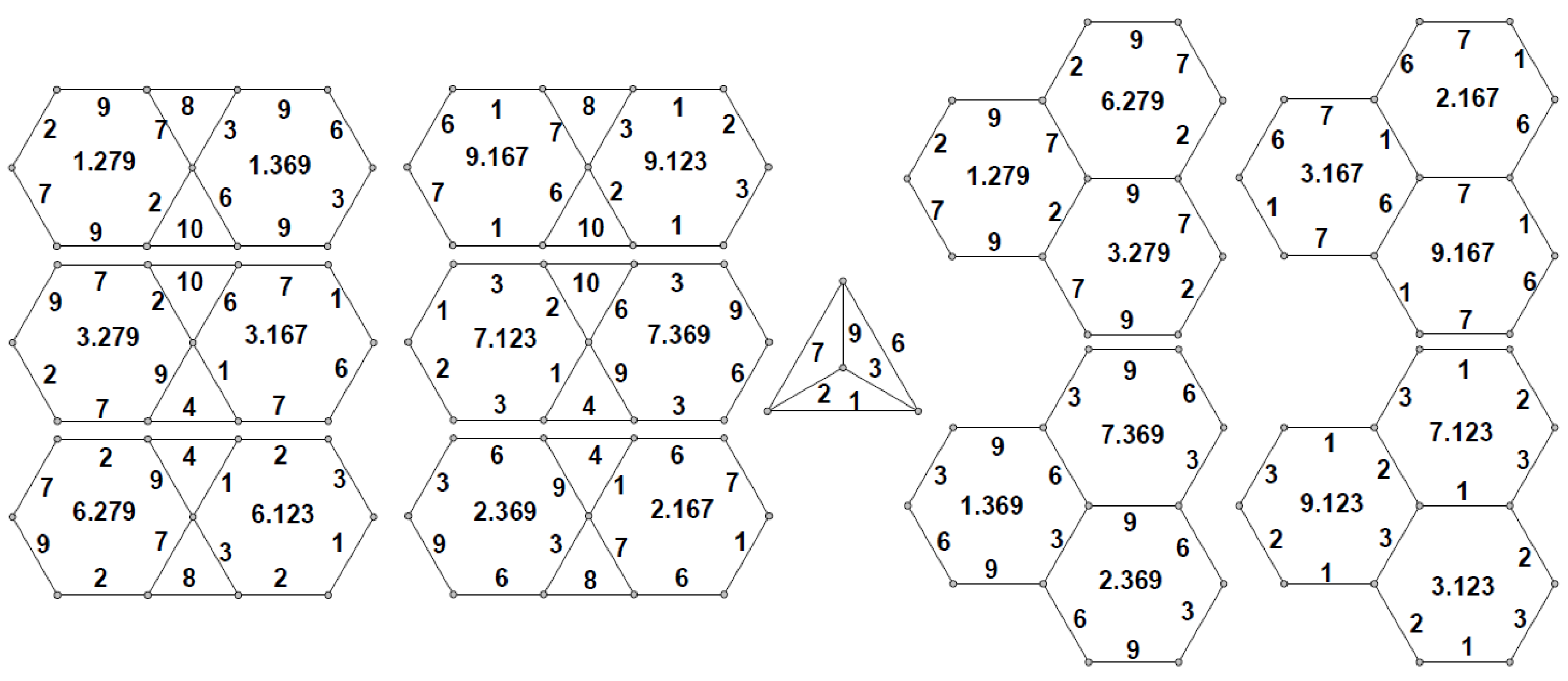}
\caption{Tessellated neighborhoods of a vertex of $G_{n,4}$ in the subgraphs of Theorem~\ref{1}}
\end{figure}

\begin{theorem}\label{1}
There exists an infinite family ${\mathcal G}_1$ of finite connected graphs $G=G_{n,4}$ with asymptotic diameter $|V(G)|^{1/3}$ such that that the subset $V_6$ of vertices $v\in V(G)$ with {\rm deg}$(v)=\Delta(G)=6$ has asymptotic order $|V(G)|$. In that case, almost every $v\in V_6$:
\begin{enumerate}\item[$1.$]
is incident to three triangles $T_0,T_1,T_2$ in $G$ with pairwise
intersection $\{v\}$ determining exactly six planar ${\mathcal D}$-modeled subgraphs $D_{i,j}^k$ $(i,j=0,1,2$; $k=0,1)$ such that $T_i\cup T_j= D_{i,j}^0\cap D_{i,j}^1$ for each pair $\{i,j\}\subset\{0,1,2\}$ with $i\ne j;$
\item[$2.$]
is the intersection of the six $\mathcal D$-modeled subgraphs of $G$ above, in which {\rm deg}$(v)=4$, and exactly four $\mathcal H$-modeled subgraphs in $G$, in which {\rm deg}$(v)=3$, and such that the closed neighborhood of $v$ in $G$ is contained in a fixed way in the union of the tessellated neighborhoods of $v$ in the ten cited subgraphs, comprising $43$ vertices.\end{enumerate}
\end{theorem}

To give an idea of what is going on locally at almost every vertex in the context of Theorem~\ref{1}, Figure 1 shows on its left (resp. right) side the closed ---tessellated--- neighborhoods of a particular vertex $v$ ---given by the edge-colored copy of $K_4$ depicted at the figure center--- in $G_{n,4}$ ---or $G_{\infty,4}$, see Section~\ref{s5}--- in each of the ten subgraphs mentioned in the two items of the statement, namely, in the six $\mathcal D$- (resp. four $\mathcal H$-) modeled subgraphs of $G_{n,4}$ claimed above, for a value of $n$ sufficiently large, with edges colored via $a=7$, $b=9$, $c=2$, $d=3$, $e=1$ and $f=6$.

\begin{cor}\label{2} There is a subfamily ${\mathcal G}'_1$
of ${\mathcal G}_1$ such that any $D_{i,j}^k$ in a member $G$ of ${\mathcal G}'_1$
is a $\mathcal D$-modeled subgraph restricted to a $30^\circ$-$60^\circ$-$90^\circ$ triangular region of the Euclidean plane. Moreover, there are $n-1$ pairwise distinct
such subgraphs $D_{i,h}^k$ distributed, for $y\ge 1$, into two subsets of size $\frac{n-1}{2}$ composed each by isomorphic subgraphs. By denoting these $\frac{n-1}{2}$-subsets by $V^-_y$ and $V^+_y$, if $k=5+2y$; resp. $U^-_y$ and $U^+_y$ if $k=4+2y$, with $|V^-_y|<|V^+_y|$ and $|U^-_y|<|U^+_y|$, then $|V^-_y|=y^2+y-1$ and $|V^+_y|=3y^2+3y-3-\epsilon(k)$, where $\epsilon(k)=1$ if $k\equiv 1$ {\rm mod 3} and $\epsilon(k)=0$ if $k\not\equiv 1$ {\rm mod 3};
resp. $|U^-_y|=|V^-_y|-y$ and $|U^-_y|=|V^+_y|-3y$.
\end{cor}

Figure 9 of \cite{De4} illustrates the $30^\circ$-$60^\circ$-$90^\circ$-triangular regions in Theorem~\ref{1}; alternatively, see Figures 6 and 7 below.
The proofs of Theorem~\ref{1} and Corollary~\ref{2} in Section~\ref{s9} are composed by
the arguments presented in Sections~\ref{s3}-\ref{s9} and, for the ${\mathcal H}$-modeled subgraphs in item 2 of Theorem~\ref{1}, by Theorem 2 of \cite{De4}.

\section{$K_3$-types and $K_3$-type graphs}\label{s3}

A triangle in $\Gamma_n$ has $K_3$-{\it type} $(a,b,c)$ if its edges
have colors $a,b,c\in[k]$. If no confusion arises, we suppress
commas and parentheses, so we write $(a,b,c)=abc$. More generally, a $K_3$-{\it type}
$abc=acb=bac=bca=cab=cba$ of $\Z_n$ is a 3-multiset $\{a,b,c\}$ of
$[k] \cup \{0\}$ such that $a+b\in\{c,-c\}\in[k]$, where $a+b$ is
taken mod $n$. (This 3-multiset can be viewed as a class of at most
six 3-tuples of colors of $[k]\cup\{0\}$, one of which is $abc$).

\begin{exm} The $K_3$-types $\{a,b,c\}$ of $\Z_7$ with $\gcd(a,b,c)=1$ are $\{0,1,1\}$, $\{1,1,2\}$, $\{1,2,3\}$, $\{1,3,-(1+3)=3\}$ and $\{2,3,-(2+3)=2\}$, where the greatest common divisor $\gcd(J)$ of a finite multiset $J$ of nonnegative integers is the largest common divisor of the nonzero integers of $J$.\end{exm}

Let $G_n$ be the graph whose vertices are the $K_3$-types of $\Z_n$
and such that any two of them, say $v$ and $v'$, are adjacent
via an edge $\epsilon$ if and only if $v$ and $v'$ share either two different
colors of $\Gamma_n$ or one color of $\Gamma_n$ repeated twice, say
$a$ and $a'$; in either case we can consider $\epsilon$ as determined
by $\{v,v'\}$ or by $\{a,a'\}$. We take $\{a,a'\}$ ($=aa'$, for
short) as the {\it color} of $\epsilon$, so that $G_n$ becomes an
edge-{\it colored} graph. In addition, we assume that $G_n$ does
not have multiple edges.
In the example above, only
123 is rainbow. Each rainbow triangle $t$ in $\Gamma_n$ and edge
$\epsilon$ of $t$ determine exactly one rainbow triangle $t' \neq t$
with the same colors of $t$ and sharing $\epsilon$ with $t$.
For $n=2k+1\geq 7$, let $G'_n\subseteq G_n$ be the subgraph of $G_n$
induced by the rainbow $K_3$-types of $\Z_n$. Let $G_{n,3}$ be the
component of $G'_n$ containing the $K_3$-type $123$. Then all the
remaining components of $G'_n$ are isomorphic to graphs $G_{m,3}$ with $1 < m < n$ and $m|n$. Notice that the
vertices of $G_{m,3}$ are 3-sets.
Now, consider $\N=\{m\in\Z:m\geq 0\}$ as an infinite {\it color
set}. A $K_3$-type $abc$ of $\Z$, simply called a $K_3${\it -type},
is a 3-multiset $\{a,b,c\}$ of $\N$ such that the sum of the two
least colors equals the greatest one. Let $G_{\infty,3}$ be the
graph whose vertices are the $K_3$-types $abc$ with $\gcd(a,b,c)=1$
and whose edges are as defined above for $G_n$.
Given $m,m',n\in\N$ with $m'\in[k]$, we say that $m'\equiv m$
MOD $n$ whenever if for $m''\equiv m$ mod $n$ with $0 \leq
m''<n$: {\bf (1)} if $m''>n/2$, then $m'=n-m''$; {\bf (2)}
if not, then $m'=m''$.
Here, $m'$ is said to be the {\it reduction of} $m$ MOD $n$.
It was shown in \cite{DHS}, Proposition 2.16, that for odd $n\geq
7$, $G_{n,3}$ can be obtained, from a connected subgraph $F$ of
$G_{\infty,3}$ containing 011, 112, 123 and the remaining $K_3$-types
with colors $\leq n$, by reducing  MOD $n$ all the colors of
$K_3$-types of $F$.
Let $\phi(n)$ be the value of Euler's totient function at the
positive integer $n$. It was shown in Theorem 2.17 of \cite{DHS}
that $|V(G_{n,3})|=O(n\phi(n))$ and subsequently, in Theorems 2.20
and 2.21, that the diameter of $G_{n,3}$ is both $\Omega(n)$ and
$O(|V(G_{n,3})|^{1/2}\,)$. The family ${\mathcal G}_0$ in the
introductory section above is formed by these graphs $G_{n,3}$.

\section{$K_4$-types and $K_4$-type graphs}\label{s4}

\begin{figure}[htp]
\hspace*{2cm}
\includegraphics[scale=0.20]{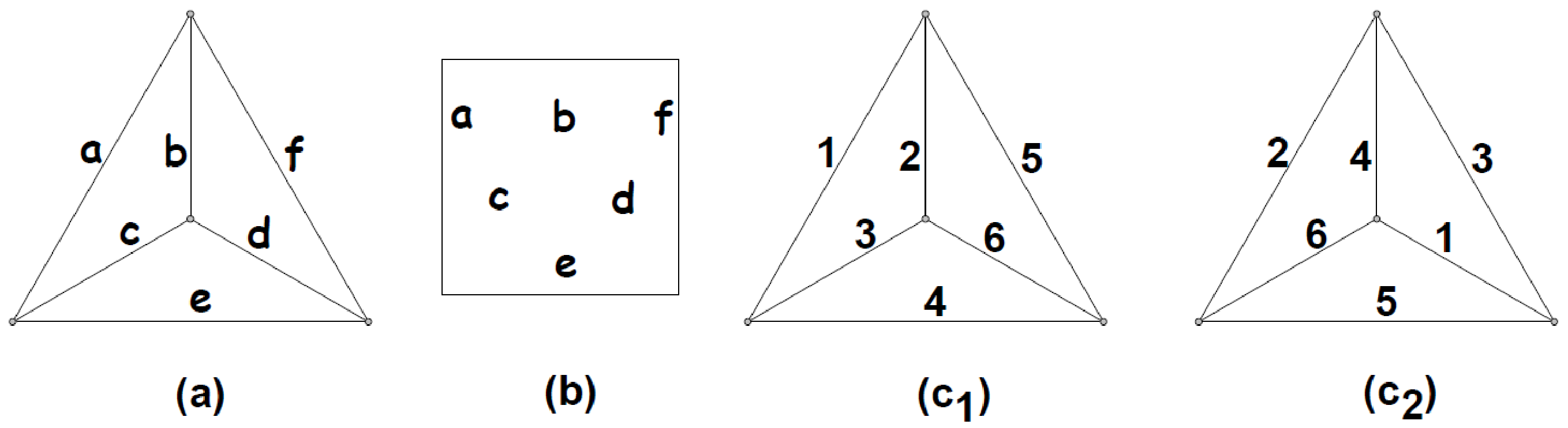}
\caption{Representing a generic $K_4$-type $abcde\!f$ and its cases MOD 13}
\end{figure}

A $K_4$-{\it type} of $\Z_n$ (resp. $\Z$) is a maximal class of
$6$-tuples $abcde\!f$ of colors of $[k]$ (resp. $\N$) such that $abc$,
$cde$, $ae\!f$ and $bdf$ are $K_3$-types of $\Z_n$ (resp. $\Z$). Such a
class has at most twenty-four 6-tuples. A 6-tuple in a $K_4$-type $t$ is
called a {\it card} of $t$. If no confusion arises, we represent a
$K_4$-type by one of its cards. A card $abcde\!f$ will be represented: {\bf (i)} either as a tetrahedron each of
whose edges bears a color, as in Figure 2(a); {\bf (ii)} or by keeping only the locations of the colors in (i) in an enclosure, as shown in Figure 2(b).

The colors in Figure 2(a) split into three different pairs of
opposite colors: $\{a,d\}$, $\{b,e\}$, $\{c,\!f\}$, (opposite in the
sense that each pair is held by a corresponding pair of edges of
$K_4$ with no vertices in common, the remaining edges forming a
4-cycle).

Any 6-multiset of $\N$ determines {\it at most} one $K_4$-type of
$\Z$. This is not true for $(\Z_n,[k])$ in place of $(\Z,\N)$. For example,
the two rainbow $K_4$-types $123645$ and $246153$ of $\Z_{13}$
represented in Figures 2(c$_1$) and 2(c$_2$), respectively, are
distinct but have the same underlying multiset.

A {\it rainbow} $K_4$-type is one with six
different colors. Given $n=2k+1\geq 13$, let $G'_{n,4}$ be the graph whose vertices
are the rainbow $K_4$-types $abcde\!f$ of $\Z_n$ with
$\gcd(a,b,c,d,e,\!f,n)=1$ and such that any two such vertices, say
$t$ and $t'$, are adjacent via an edge $\epsilon$ if and only if $t$
and $t'$ looked upon as $K_4$-types share precisely two $K_3$-types
$v$ and $v'$. In this case, $v$ and $v'$ share exactly one color $a$
of $[k]$. We take $a$ as the ({\it weak}) {\it color} of $\epsilon$ and this makes
$G'_{n,4}$ into an edge-colored graph.

In order to distinguish the $\mathcal D$- and $\mathcal H$-modeled subgraphs
that we claim  $G'_{n,4}$ contains, we introduce
the graph $G''_{\infty,4}$ as the simple graph (i.e., graph without loops or multiple edges) whose
vertices are the $K_4$-types $abcde\!f$ with $a\neq d$, $b\neq e$
and $c\neq f$ unless $abcde\!f=011011$ and satisfying
$\gcd(a,b,c,d,e,\!f)=1$, with two vertices $u$ and $v$ determining
an edge if and only if they share precisely two $K_3$-types in differing locations of the representation of the $K_4$-types that stand for $u$ and $v$ as in Figure 1.

\begin{figure}[htp]
\hspace*{2cm}
\includegraphics[scale=0.2]{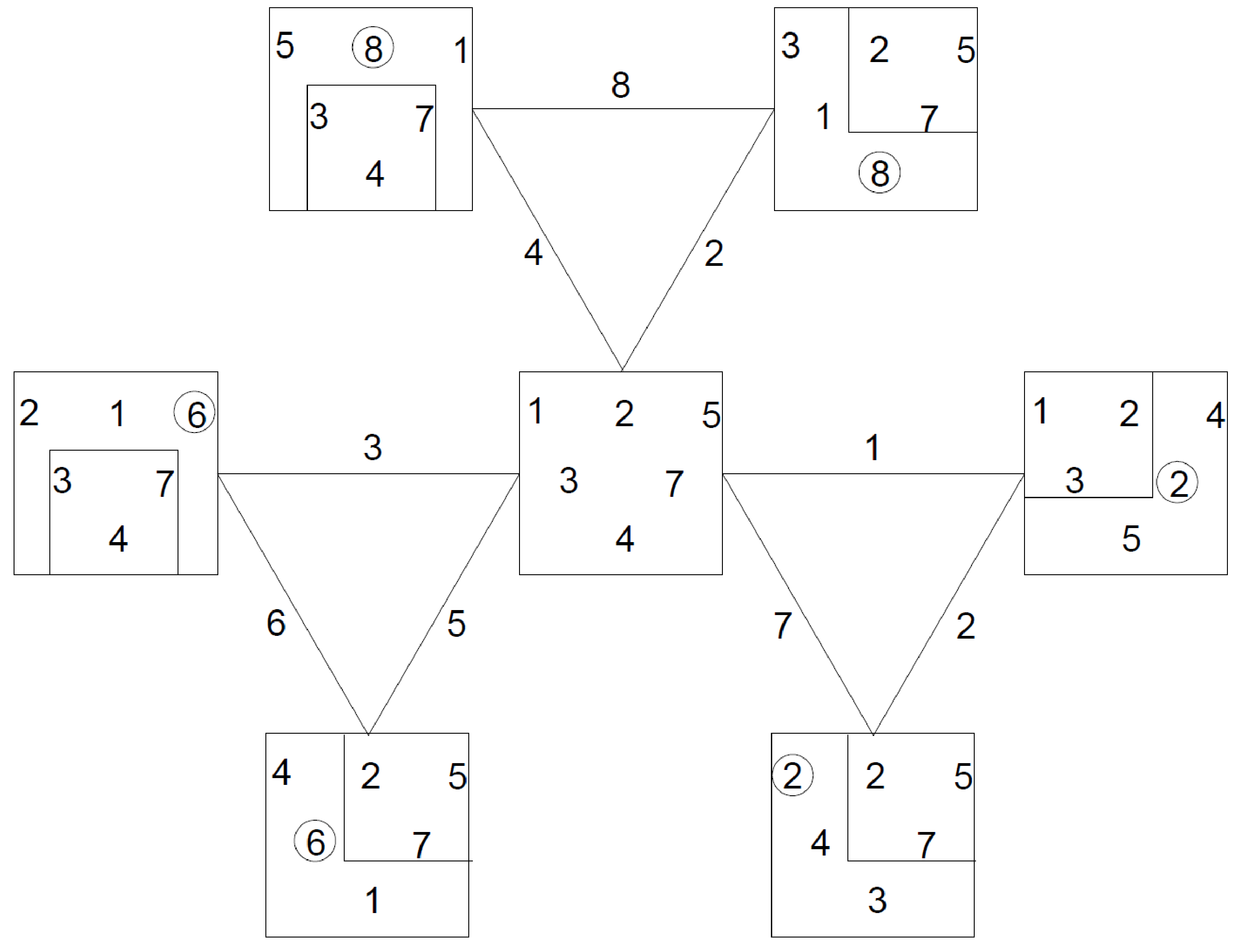}
\caption{A neighborhood of $123745$ in $G''_{\infty,4}$}
\end{figure}

Figure 3 illustrates $G''_{\infty,4}$ as well as Theorem~\ref{t} below. The figure
represents a neighborhood $N$ of the $K_4$-type $123745$ in
$G''_{\infty,4}$. Notice that the two right-lower
$K_4$-types in Figure 3 (joined by the edge colored with 6) are not rainbow. An edge $\epsilon$ joining two
vertices $t$ and $t'$ of $G''_{\infty,4}$ with respective cards $r$
and $r'$ determines a $K_3$-type $s$ common to $t$ and $t'$ and {\it
equally located} in $r$ and $r'$ in the sense that the component
colors of $s$ occupy the same positions in $r$ and $r'$, just as the
$K_3$-type $s=123$ is not only common to but also equally located in
the central card in Figure 3 and the card horizontally located at
its right, with $s$ occupying the three upper-left locations
in $r$ and $r'$. The {\it locations} $g_r$ of
the colors in the cards $r'$ of the statement of Theorem~\ref{t}
obtained from the central card $r$ at the center of Figure 3 are
shown encircled. Also, the $K_3$-type $s$ is highlighted in a
sub-enclosure of its own. Observe that in each of the six
enclosures representing the neighbors of the central vertex in
Figure 3 the two colors outside the sub-enclosure and the encircled
color are permuted in their positions.

\begin{theorem}\label{t}
Let $t\in V(G''_{\infty,4})$. Let $r$ be a card of $t$ with color
$g$ at location $g_r$ and color $g'$ at the location $g'_r$ opposite
to $g_r$. Then $t$ has a neighbor $t'$ with card $r'$ differing from
$r$ just in: {\bf (a)} the color at $g_r$ and {\bf (b)} a
permutation of the colors at the two locations $\neq g'_r$ in just
one of the two $K_3$-types common to $r$ and $r'$ that contain the
color at $g_r$.\end{theorem}

\begin{proof} $t'$ is determined from $t$ as follows. Let $s,s'$ be the
two $K_3$-types not containing $g_r$ in $r$. Then $s$ and $s'$ contain
$g'_r$. We can assume that $s'$ has its colors equally located in $r$
and $r'$. Let $i,j$ be the colors of $r$ at the two locations
$i_r\neq g'_r$ and $j_r\neq g'_r$ of $s$. Thus $s=ijg'$. The two other
$K_3$-types in $t$ apart from $s$ and $s'$ are of the form $gij'$
and $gji'$ with $s'=i'j'k$.  We take $r'$ as having the colors
$i,j$ exchanged with respect to $r$. So
$(i_{r'},j_{r'})=(j_r,i_r)$. Let $\nu(a,b) = \{\mid a - b \mid \}
\cup \{ a + b \}$ for each pair of integers $a,b\geq 0$. There is
at least one color $h\in\nu(i,j) \cap \nu(i',j')\neq\emptyset$
that yields $r'$ when located at $g_r$ (which should be called
$h_{r'}$ in $r'$) so that $r'$ is formed by the $K_3$-types
$s=ijg'$, $s'=i'j'g'$, $hii'$ and $hjj'$. Moreover, $r'$ does not
depend on the selected card $r$ of $t$. In fact $h=h(r,g_r)$
depends only on $r$ and $g_r$. If $r=011011$ and $g=0$ then $h$
equals either 0, yielding $t'=t$, not a distinct neighbor of $t$
in $G''_{\infty,4}$ so we discard it, or 2, yielding a neighbor
$t'$ of $t$. Otherwise, since no remaining vertex of $G''_{\infty,4}$
is of the form $abcabc\neq 011011$, then $|\nu(i,j) \cap \nu(i',j')|=1$, even if
$(r,g)=(011011,1)$. Thus, if either $r\neq 011011$ or
$(r,g)=(011011,1)$, then $h$ is unique. \end{proof}

\begin{exm} In the following special cases, $g$ assumes
subsequently colors $f$, $a$ and $d$ in a $K_4$-type $t$ of card
$r=abcde\!f$: {\bf (A)} applying Theorem~\ref{t} to
$(r,g)=(112354,4)$ (so $g=f$) yields $t'=t$ where $g_r=f_r=4_r$
because exchanging $d_r=1_r$ and $e_r=1_r$ does not produce changes
from $r$; {\bf (B)} applying Theorem~\ref{t} to $(r,g)=(011011,0)$
(so $g=a,d$) yields, for $h=2$, neighbors $t',t''$ with respective
cards $r'=211011$ and $r''=011211$ where $g_r=a_r,d_r$
respectively, but observe that $t'=t''$. \end{exm}

\section{Canonical triangles}\label{s5}

Let $G_{\infty,4}$ be the supergraph of $G''_{\infty,4}$ obtained by adding to the vertices of $G''_{\infty,4}\backslash\{011011\}$ the loops offered by the method of vertex adjacency in Theorem~\ref{t} and Figure 3, taking each maximal set of loops incident to a common vertex and with a common color to have multiplicity
1. Then, a link or loop joining vertices $t$ and $t'$ in
$G_{\infty,4}$ has the pair $(s,s')$ in the proof of Theorem~\ref{t} as its {\it strong color} and the only color $g'$ in $s$ and $s'$ that remains at
the location $g'_r=g'_{r'}$ both in $r$ and $r'$ as its {\it weak color}. Let
$G'_{\infty,4}$ be the graph obtained from $G_{\infty,4}$ by restriction to the vertices that are rainbow $K_4$-types.

Applying {\rm Theorem~\ref{t}} to the colors $g,g'$ of a pair of opposite edges of a vertex $t$ of $G_{\infty,4}$ looked upon as a $K_4$-type with card $r$ yields $h(r,g)=h(r,g')$. This determines in $r$ two corresponding neighboring cards $r'$ and $r''$ representing respective neighbors $t'$ and $t''$ of $t$. The two $K_3$-types that $r'$ and $r''$ share and those two that $r$ and $r'$ (resp. $r$ and $r''$) share constitute the four $K_3$-types of $r'$ (resp. $r''$).
The resulting triangle, whose vertices $t,t',t''$ have respective cards $r,r',r''$, is said to be a {\it canonical triangle}, or {\it CT}. Since there are three pairs of opposite vertices in the card $r$ associated to the vertex $t$ of $G_{\infty,4}$, then there are at most three CTs incident to $t$. Since each $G'_{n,4}$ can be obtained from $G'_{\infty,4}$ via reduction MOD $n$, we have completed the proof of the following corollary.

\begin{cor}\label{xx}
The graphs $G'_{\infty,4}$ and $G'_{n,4}$ are edge-disjoint
unions of {\rm CT}s, at most three such {\rm CT}s incident to each vertex.\hfill\rule{2mm}{2mm}
\end{cor}

When two or three $K_4$-types in a CT $T=\{t,t',t''\}$ obtained as in Theorem~\ref{t} coincide (e.g., either
$t=t'\neq t''$ or $t=t''\neq t'$ or $t\neq t'=t''$ or $t=t'=t''$),
then we say that $T$ is a {\it degenerate CT}.

\begin{exm} {\bf(A)} If $t$ has $r=abcde\!f$ with $a,b>0$,
$c=a+b$, $d=a$, $e=b$, $f=|a-b|$ and
$(g_r,g'_r)\in\{(a_r,d_r),(b_r,e_r)\}$, then $t'=t''$. This
yields two degenerate CTs with vertices of the form $t$, $t'$ and $t''=t'$, where $tt'=tt''$ and $t't''$ is a loop of $G_{\infty,4}$.
{\bf(B)} Theorem~\ref{t} applied to $t=000111$ yields three
degenerate CTs, each representable by: two vertices, namely
$t$ (twice) and $t'=011011$, a link $tt'$ and a loop at $t$; these
three CTs coincide, since edges are assumed to have multiplicity 1.
{\bf(C)} Theorem~\ref{t} applied to $t=132112$ yields three CTs incident to $t$, one of which,
obtained by making value changes in both cases of color $g=2$ at opposite locations in $t$, has its three vertices equal to $t$, so this CT reduces to a looped vertex in $G_{\infty,4}$. The two remaining CTs incident to $t$ are $\{t,202111, 132201\}$ and $\{t,431122, 132421\}$.
\end{exm}

\begin{cor}\label{conn} $G_{\infty,4}$ is connected.\end{cor}

\begin{proof} Given $t=abcde\!f$ and $t'=abcydx$ in $G_{\infty,4}$ there exists a 2-path
in $G_{\infty,4}$ from $t$ to $t'$ with   middle vertex card
$abc\!f\!xd$ and edge strong colors $\{abc,bdf\}$ and $\{abc,adx\}$.
Let $cde$ and $cxy$ be $K_3$-types of $\Z$ with
$\gcd(c,d,e)=\gcd(c,x,y).$ Then there exists a path in
$G_{\infty,4}$ whose ends have cards of the form $abcde\!f$ and
$abcxyz$. This uses the fact that if $\gcd(c,d,e)=\gcd(c,x,y)$, then
there is a path in $G_{\infty,3}$ from $cde$ to $cxy$ \cite{DHS}.
Thus, if $abcde\!f\in V(G_{\infty,4})$, then there exist: {\bf (a)}
a path in $G_{\infty,4}$ from $110110$ to $110aa(a+1)$; {\bf (b)} a
path in $G_{\infty,4}$ from $110aa(a+1)$ to $aa0bbc$; {\bf (c)} a
path in $G_{\infty,4}$ from $aa0bbc$ to $abcde\!f$. Hence, every
vertex of $G_{\infty,4}$ can be connected to $110110$. \end{proof}

\section{Generation of ${\mathcal D}$-modeled subgraphs}\label{s6}

\begin{cor}\label{tt}
The set of {\rm CT}s of $G_{\infty,4}$ is in $1$-$1$ correspondence with the
family of $4$-multisets or quadruples $abcd$ of colors of $\N$ such
that: {\bf (a)} $\nu(a,b) \cap \nu(c,d) \neq \emptyset$ {\rm(}or
$\nu(a,c) \cap \nu(b,d) \neq \emptyset$ or $\nu(a,d) \cap \nu(b,c)
\neq \emptyset${\rm);} {\bf (b)} $\gcd(a,b,c,d) = 1$, so at least one
of $a,b,c,d$ is nonzero.\end{cor}

\begin{proof} From Theorem~\ref{t} and Corollary~\ref{xx}, each CT of
$G_{\infty,4}$ has its vertices as $K_4$-types sharing precisely
four colors as in the statement.\end{proof}

\begin{exm}\label{11} In Figure 3, the upper (resp. lower-left, lower-right) CT has
its vertices sharing the quadruple 1357 (resp. 1247, 2345).
\end{exm}

From now on, each CT will be denoted by its associated multiset
in Corollary~\ref{tt}. Given a rainbow $K_4$-type
$t=abcde\!f$, the CTs incident to $t$ are obtained by deleting
from $t$ each one of the three pairs $ad$, $be$ and $cf$, which yields respectively
$bce\!f$, $acdf$ and $abde$.

Let $abcde\!f$ be a vertex of $G_{\infty,4}$ and let $C=acdf$ and $D=abde$ be two CTs in $G_{\infty,4}$ sharing just $abcde\!f$. Then $C\cup D$ is represented as a colored 5-vertex plane graph $B(t,a,d)$ where C and D participate as respective equilateral triangles
$\x{C}$ and $\x{D}$, respectively, that share solely a vertex $t$ (i.e., $\x{C}\cap\x{D}=\{t\}$) that stands for $abcde\!f$ and is center of a point symmetry that takes $\x{C}$ onto $\x{D}$ and viceversa. Thus, pairs of sides of $\x{C}$ and $\x{D}$ incident to $t$ are set collinearly as in Figure 4. We require $a$ to tag the centers of both $\x{C}$ and $\x{D}$, and the remaining colors of $C$ and $D$ to tag respectively the vertices of $\x{C}$ and $\x{D}$ internally. Then,
$d$ is the color tagging $t$ internally in both $\x{C}$ and
$\x{D}$. We tag each edge of $\x{C}$ (resp. $\x{D}$) with the weak color of the corresponding edge of $C$ (resp. $D$), such that the weak color of each edge $\epsilon$ of $\x{C}$ forms: {\bf (a)} a $K_3$-type $s(\epsilon)$ with the colors tagging the endvertices of $\epsilon$ in $\x{C}$; {\bf (b)} another $K_3$-type $s'(\epsilon)$, with the central tagging color of $\x{C}$ and the color tagging the vertex opposite to $\epsilon$ in $\x{C}$. Notice that $\{s(\epsilon),s'(\epsilon)\}$ is the strong color of the image of $\epsilon$ in $G_{\infty,4}$. Let $\epsilon_C$ and $\epsilon_D$ be edges of $\x{C}$ and $\x{D}$, respectively, meeting at an angle of $120^\circ$ at vertex $t$. Then the color $d$ tagging $t$ in both $\x{C}$ and $\x{D}$ forms with the colors tagging $\epsilon_C$ and $\epsilon_D$ the $K_3$-type $s(\epsilon_C)=s(\epsilon_D)$.

\subsection{Growth of a $\mathcal D$-modeled subgraph}\label{pr}

\begin{figure}[htp]
\hspace*{.1mm}
\includegraphics[scale=0.29]{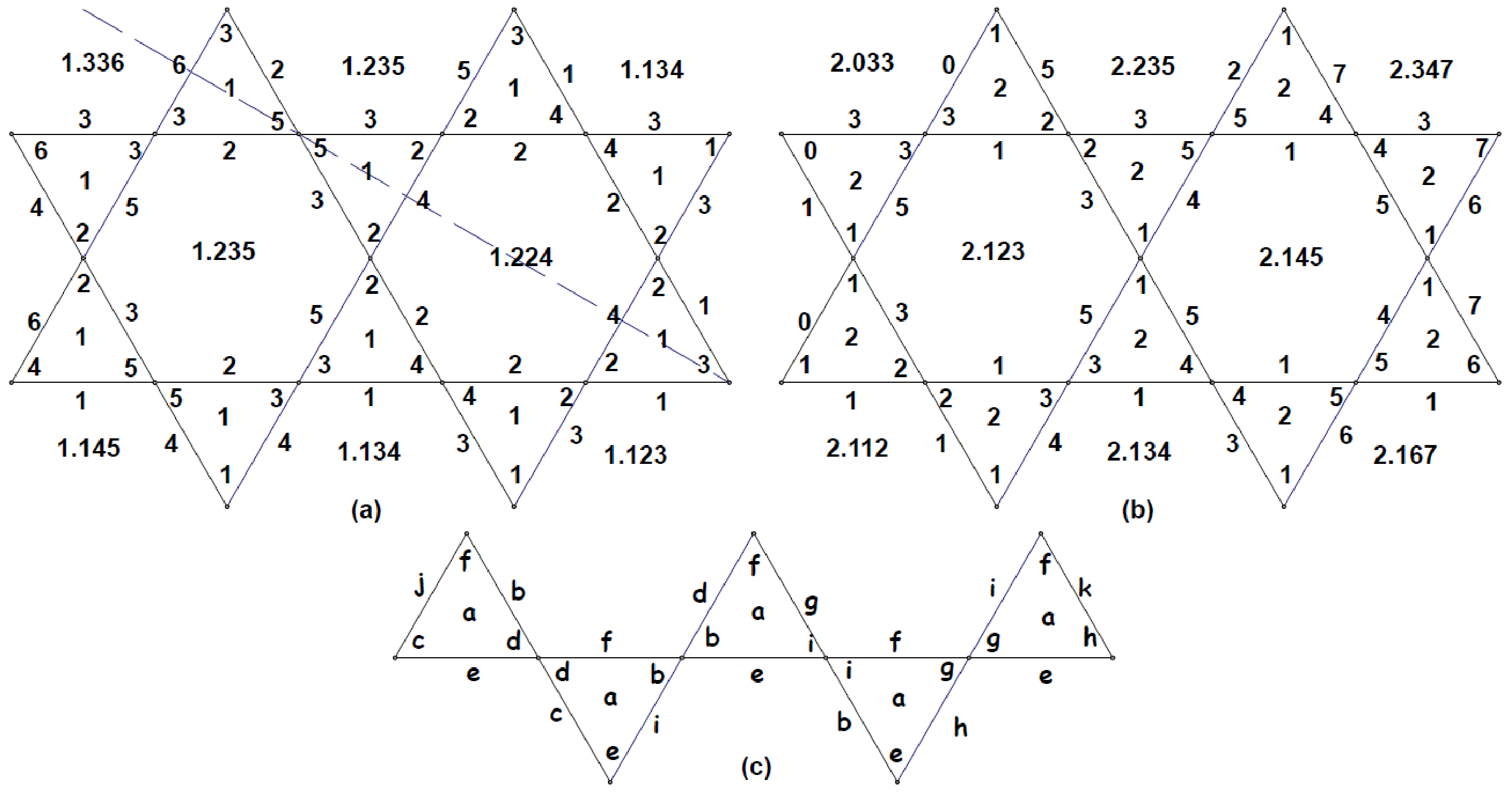}
\caption{Unfoldings of subgraphs of $G_{\infty,4}$}
\end{figure}

The growth of a $\mathcal D$-modeled subgraph of $G_{\infty,4}$ sprouting from $B(t,a,d)=\x{C}\cup\x{D}$ via Theorem~\ref{t} can be performed via the following properties deducible via Theorem~\ref{t} and enjoyed by the objects conceived in the previous paragraph with their tagging notation around $r=abcde\!f$ as shown in Figure 4(c) and  illustrated in Figure 4(a)-(b).

{\bf (1)} Given a CT $C=a\!f\!gh$, let $a$ be the
central tag of $\x{C}$ and let color $f$ tag a
vertex $u$ in $\x{C}$. Then there is a color $i$ so that {\bf (a)}
$\nu(a,h)\cap\nu(f,g)=\{i\}$; {\bf (b)} the edges $\epsilon=uu'$ in
$\x{C}$ with $u'$ having tag $g$ or $h$ in $\x{C}$ have color $i$, denoted $\gamma(\epsilon)=i$. {\bf (2)} Let $\ell$ be the line containing
$u$ and parallel to the unique edge of $\x{C}\backslash u$. Then each
pair $(u,C)$ determines at most one remaining CT $D\neq C$ sharing $u$
with $C$, so that $\x{D}=\rho_\ell(\x{C})$, where
$\rho_\ell$ is reflection of the plane on $\ell$, and
having {\bf (a)} $a$ as central tag; {\bf (b)} the tag $f$ of
$u$ in $\x{C}$ as tag of $u$ in $\x{D}$; {\bf (c)}
for each edge $\epsilon=uu'$ of $\x{C}$:
{\bf (i)} $\gamma(\epsilon)$ as the tag of $\rho_\ell(u')$ in
$\x{D}$ and {\bf (ii)} the tag of $u'$ in $\x{C}$
as the tag of $\rho_\ell(\epsilon)$. {\bf (3)} The vertex $u$ is the
$K_4$-type formed by the $K_3$-types determined by each
edge $\epsilon$ of $\x{D}$ incident to $u$ and formed by:
{\bf (a)} $a$ and the tags of $\epsilon$ and the vertex
opposite to $\epsilon$ in $\x{D}$; {\bf (b)} the tags
of $\epsilon$ and the endvertices of $\epsilon$ in
$\x{D}$.

The union of two CTs $C$ and $D$ that share precisely one vertex $v$
is said to be a {\it butterfly} and denoted $CvD$. In this case,
$v$ is called the {\it central vertex} of $CvD$. Note that the
colors of $v$ in $\x{C}$ and $\x{D}$ equal a fixed
color $d$ which we call the {\it butterfly color} of $CvD$. For
example, $B(t,a,d)$ above is a butterfly $CtD$ with central color
$a$ and butterfly color $d$, say with $C=acdf$ and $D=abde$.
Given a simple graph $G$ and a pseudograph $H$ (i.e., $H$ is a non-simple graph in which each vertex may be incident to one or more loops),
then $G$ is an {\it unfolding} of $H$ if there exists a surjective map $f:V(G)\rightarrow V(H)$ such that for each $v\in V(G)$ there exists a 1-1 correspondence induced by $f$ from the links incident to $v$ in $G$ to the edges incident to $f(v)$ in $H$.

\subsection{Maximal $\mathcal D$-modeled graphs}\label{tt1}

Let\, $t=abcde\!f$ be a rainbow $K_4$-type. A maximal $\mathcal D$-modeled graph
$H'=H'(t,a)=H'(t,a,d)\supset B(t,a,d)$ that is an unfolding of an edge-disjoint
union $H=H(t,a)=H(t,a,d)$ of butterflies in $G_{\infty,4}$ with common
central color $a$ is generated by repeated application of item $(2)$, Subsection~\ref{pr}, at gradients $0^\circ$, $60^\circ$, $120^\circ$, $180^\circ$, $240^\circ$, $300^\circ$ of the line $\ell$ in the item.

\begin{exm} Both Figure 4(a) and 4(b) show parts of an $H'$ as above.\end{exm}

We will see that if such an $H'$ is not a
subgraph of $G_{\infty,4}$, then it can be {\it folded} along at most two
{\it symmetry axes}, or {\it SA}s, to yield $H$. The dotted line in Figure
4(a) represents such an SA. In particular, edge colors will coincide by
reflection in an SA. The graph obtained from $H$ by removing the resulting loops will be seen to be a subgraph of $H'$
spanning a connected region of the plane delimited by SAs.
Edges crossing an SA at $90^\circ$ will yield loops of $H$ and each
CT in $H'$ will be incident to three hexagons.

\begin{obs}~\label{tos}
Given a vertex $t$ of $H'(t,a,d)$, the three {\rm CT}s incident to $t$
according to {\rm Theorem~\ref{t}} are: {\bf (a)} the two {\rm CT}s incident
to $t$ in $H'(t,a,d)$ and {\bf (b)} the {\rm CT} formed by the colors of
the four edges of the two {\rm CT}s in item {\rm(a)} which are incident to
$t$.\end{obs}

\section{Presence and properties of $6$-cycles}\label{s7}

The graph $H'(t,a,d)$ in Subsection~\ref{tt1} has two edge-disjoint 6-cycles with just the vertex $t$ in common which are given by regular hexagons in the plane when the CTs of $H'(t,a,d)$ are represented as equilateral triangles as in the discussion after Example~\ref{11}.
This is the specific case in Subsection~\ref{he}
below. If $q$
is any of these 6-cycles, then its edges are colored with the
component colors of a $K_3$-type $s$. In that case, we denote
$q=a.s$, where $a$ is the central color of the six CTs adjacent to
$q$.

\subsection{A procedure to determine $6$-cycles}

Let $bdf=s$ and $cde=s'$ be $K_3$-types, where $t=abcde\!f$ is a vertex of $H'(t,a,d)$.
We will see that there exists a $6$-cycle
$(t^0,t^1,t^2,t^3,t^4,t^5)$ in $H'(t,a,d)$ containing $t=t^0$.
It will be determined by the following procedure that
yields $t^i$ when $t^{i-1}$ is given, successively for $i=1,2,3,4,5$, (and
returns to $t^0=t^i$ from $t^5=t^{i-1}$, if $i=6\equiv 0$ with
indices taken mod 6):

{\bf (a)} Declare the card $r^i$ of the
$K_4$-type $t^i$  to have color $a$ (as in Figure 2(b)) fixed in the
location $a_{r^0}$ (so that $a_{r^i}=a_{r^0}$) during the entire
procedure;
{\bf (b)} denote locations $b_{r^i}=b_{r^0}$,
$c_{r^i}=c_{r^0}$ and $e_{r^i}=e_{r^0}$ regardless of changes in
their color values from the initial ones, namely $b$, $c$ and
$e$ respectively along the running of the procedure;
{\bf (c)}
define color $h^i=b$ (resp. $h^i=f$) if $i$ is even (resp. odd);
{\bf (d)}
establish a color exchange via a redesignation of
locations at the $i$-th level:
$d_{r^i}=h^{i-1}_{r^{i-1}}\,\,\,\mbox{  and
}\,\,\,h^i_{r^i}=d_{r^{i-1}};$
{\bf (e)} the color $e_{r^i}$
(resp. $c_{r^i}$) if $i$ is even (resp. odd) takes the only value from
$\nu(a_{r^i},\!f_{r^i})\cap\nu(c_{r^i},d_{r^i})$
(resp. $\nu(a_{r^i},b_{r^i})\cap\nu(d_{r^i},e_{r^i})$.) This determines a
well-defined card $r^i$ and yields a location instance for
the determination of a 6-cycle as claimed.

\begin{exm} A 6-cycle generated by the procedure in the previous paragraph and starting at $t^0=123745$ is
$$a.s=1.257=(123745, 123587, 156287, 156712, 176512, 176245).$$ Its
accompanying coplanar 6-cycle $a.s'$ is
$$1.347=(123745,187345,187434,134734,134376,123476).$$ An
essentially equivalent 6-cycle to this and sharing its first two
vertices with $a.s'$ as just given is $7.145=(123745, 583741,
48C751, 1BC754, 5B6714, 426715)$, where capital
hexadecimal notation is used, and its accompanying coplanar
6-cycle is $7.123=$ $$(123745,321785,23178A,13279A,312796,213746),$$ sharing its first two vertices with $a.s$.
\end{exm}

\subsection{On $6$-cycles containing specific $K_4$-types}\label{he}

Each $t$ as above is
contained in precisely two $6$-cycles $q=a.s$ and $q'=a.s'$ of $H'(t,a,d)$. The
edge-color sets of $q$ and $q'$ are respectively $\{b,d,\!f\}$ and
$\{c,d,e\}$, each color tagging opposite edges.
Moreover, the color tagging $t$ in its incident CTs in $H'(t,a,d)$
and those tagging the two edges in $q$ $($resp. $q')$ that are incident to $t$
conform $s$ $($resp. $s')$. Furthermore,
$d$ is the color tagging $t$ in its incident CTs in $H'(t,a,d)$ as well as tagging
the two parallel edges of $a.bdf$ $($resp. $a.cde)$ incident neither to $t$
nor to its corresponding opposite vertex.

Given $K_3$-types $bcd$ and $bc'd'$ with $b<c<d$ and $b<c'<d'$,
define $bcd<bc'd'$ if and only if $c+d<c'+d'$. A graph $H'=H'(t,a,d)$ as in
Subsection~\ref{tt1} is said to be a {\it T-subgraph} and denoted $a(s)$,
where $s$ is the smallest $K_3$-type $\neq 000$ coloring a
$6$-cycle of $H'$ under '$<$', while $H=H(t,a,d)$ is denoted
$a[s]$. Hexagons $a.s$ of an $H'(t,s,d)$ and their images in
$H(t,a,d)$ are called {\it canonical hexagons} or {\it CH}s.

\begin{prop}\label{dorm}
Let $H'=H'(t,a,d)$, where $t=abcde\!f$ is common to $C=acdf$ and
$D=abde$, with $\x{C}\cup\x{D}\subset H'(t,a,d)$ and
$d$ tagging $t$ in both $\x{C}$ and $\x{D}$. Then, the
T-subgraph $H''=H'(t,d,a)$ has $t$ common to a
flipped copy $\x{\x{D}}$ of $\x{D}$ and a
direct copy $\x{\x{C}}$ of $\x{C}$. As
a result, $d.ca\!f$ and $d.bae$ contain the colors of the {\rm CT}s incident
to $t$ in $H''$. Moreover, $H''=H'$ if and only if $f=c$ and
$e=b$.
\end{prop}

\begin{figure}[htp]
\hspace*{.5cm}
\includegraphics[scale=0.30]{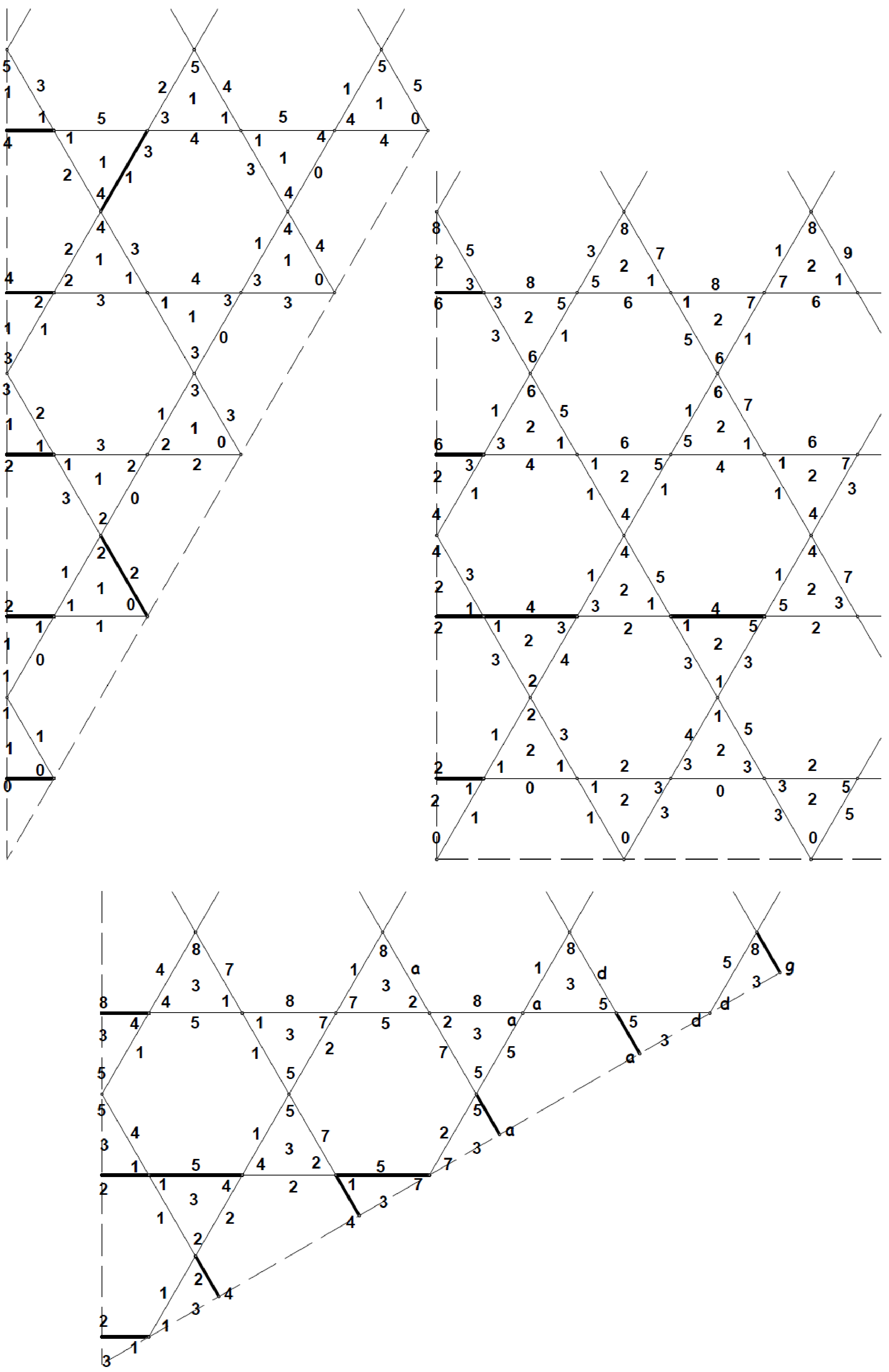}
\caption{Charts for $1[011]$, $2[011]$ and $3[112]$}
\end{figure}

\begin{proof} $H''=H'(t,d,a)$ is established as follows: {\bf (1)} represent
$H''$ as a temporarily uncolored T-subgraph and set $t$ as one of
its vertices; {\bf (2)} represent $\x{\x{C}}$ and $\x{\x{D}}$ in
$H''$ as the respective CTs $\x{C}$ and $\x{D}$ of $H'$ with common
vertex $t$ but set the locations of $a$ and $d$ in $\x{\x{C}}$ and
$\x{\x{D}}$, instead, as those of $d$ and $a$ in $\x{C}$ and
$\x{D}$, respectively; {\bf (3)} the vertex colors $c$ and $f$ in
$\x{\x{C}}$ are exchanged with respect to their locations in $\x{C}$
while the two vertex colors $b$ and $e$ in $\x{\x{D}}$ are left as
in $\x{D}$. The remaining colors of $H''$ can be set uniquely as in
Subsection~\ref{pr} above. If $H''\neq H'$, then reflection with respect
to the line perpendicular to the line $\ell$ in Subsection~\ref{pr}
through $t$ takes each edge color of $\x{\x{D}}$ in $H''$ to its
location in $\x{D}$, while the edge colors of $\x{\x{C}}$ remain as
in $\x{C}$. The statement follows immediately, as illustrated in
Figure 4, where (b), at right, represents part of the T-subgraph
$H''$ corresponding to the T-subgraph $H'$, partly represented
itself in (a), with $t=235142$ at the center in both
representations. \end{proof}

\section{From $\mathcal D$-modeled subgraphs to charts}\label{s8}

Local plane representations of some subgraphs $a[s]=a[bcd]$ of $G_{\infty,4}$ are provided in Figure 5 with notation given before Proposition~\ref{dorm}, $a=10,d=13,g=16$ and thin (resp. thick) edges for links (resp. loops). In fact, the subgraphs induced by the set of links of these $a[s]$ yield subgraphs of
the corresponding graphs $a(s)=a(bcd)$. Concretely, Figure 5
upper-left (resp. upper-right) shows a plane region delimited by two
dotted lines $\ell$ and $\ell'$ that form an internal angle of $30^\circ$
(resp. $90^\circ$) and determine a partial representation of
$H'(s,1)=1(011)$ (resp. $H'(s,2)=2(011)$), where $s=110001$
(resp. $s=211011$). This representation can be identified with
$H(s,1)=1[011]$ (resp. $H(s,2)=2[011]$) by interpreting as a loop each
thick edge interrupted perpendicularly by some dotted line
$\ell$. Moreover, $H'(s,1)$
(resp. $H'(s,2)$) is obtained by unfolding $H(s,1)$ (resp. $H(s,2)$) along
the SAs formed by the lines in the finite sequence $\ell_0=\ell$,
$\ell_1=\ell'$, $\ldots$, $\ell_i=$ reflected line of $\ell_{i-2}$
on the line $\ell_{i-1}$, for $i=2,\ldots, k-1$, where additionally
$\ell_{k-1}=$ reflected line of $\ell_1$ on the line $\ell_0$,
with $k=360^\circ/30^\circ=12$ (resp. $k=360^\circ/90^\circ=4$).

The extensions of these partial pictures to the plane will be referred to as {\it charts}. Observe that the two charts in the previous paragraph are the only charts of the form $H'(t,a)$ with $a=1,2$. However, no remaining value of $a$ produces just one chart. For example, there are two charts
$H'(s,3)$, one of which is $3(112)$, with $3[112]$ partially shown
in the bottom of Figure 5, where two straight lines $\ell_0$ and
$\ell_1$ at an angle of $60^\circ$ delimit its representation, and with
finite sequence $\ell_0,\ell_1,\ldots$, as above, of length $k=360^\circ/60^\circ=6$. The remaining
$H'(s,3)$ is $3(011)$, with $3[011]$ having exactly one SA,
delimiting a semi-plane representation.
As $a$ increases its value, the first chart $H$ not having an SA is $H=6(123)=6[123]$.

\subsection{Unfolding charts}\label{sym}

To see how the unfolding of a graph $a(bcd)$ onto its
corresponding $a[bcd]$ takes place, we observe that
if $H(t,a)\neq H'(t,a)$, then $H(t,a)$ is obtained by folds of
$H'(t,a)$ along SAs of two types:\begin{enumerate}
\item SAs dividing all CHs
of the form $a.0cc$ in symmetric halves through vertices colored
with $0$ in CTs of the form $a0cd$, i.e., through all vertices of
the forms $0bbcca$ and $0ccdda$;
\item SAs dividing all CHs of
the form $a.bbc$ in symmetric halves and passing at $90^\circ$
through the midpoints of their edges colored with $c$
(which are thick edges that yield loops) and through the
vertices opposite to them in corresponding CTs.\end{enumerate}
Here, only the CT of the form $3(123)$ has two such SAs.

In a chart $H'$, a thick edge halved
perpendicularly in its middle point by some SA yields a {\it
half-edge} of $H$, and a CT that contains a half-edge yields
a {\it half-CT} of $H$. Degenerate CT 1113, shown in the
lower-left corner of the chart $3[112]$ in Figure 5, has its center as the
intersection of two SAs (and three SAs in $3(112)$) and constitutes the
only {\it one-sixth-CT} of any chart of $G_{\infty,4}$. See also the example (C) before Corollary~\ref{conn} in Section 5, where the CTs in their shown order are
1113, 1122 and 1123, the first two present in $3[112]$. The following properties are observed:
\begin{enumerate}
\item A maximal connected region of an $H'(t,a)$ delimited by
SAs but with its interior not intersecting any remaining SA yields a
chart of $H(t,a)$.
\item Charts $a(bcd)$ and $a[bcd]$ exist, for
$b\leq c\leq d$, if and only if $c+d\leq a$.
\item Every loop of $G_{\infty,4}$ not in CTs $0011, 1111, 0112, 1113$ appears as a
half-edge in two different charts and as a thick edge in a different one. The CT that
contains such a loop:
{\bf (a)} is of the form $aabc$, where $a,b,c$ are pairwise different and $(2a,b,c)$ is a
$K_3$-type; {\bf (b)} appears as a half-CT obtained by
halving a degenerate CT as in the example (A) in Section 5 by means of an SA
in $b[112]$ or $c[112]$, and as a $3$-cycle in $a[011]$.
\end{enumerate}
Two edges in a butterfly $B(t,a,d)$ are said to
be {\it opposite} if none has $t$ as an
endvertex. Each butterfly has just one pair of opposite
edges.

\subsection{Color-alternating infinite paths}\label{17}

Any infinite path of $H'=H'(t,a)=a(bcd)$ contained in a line has successive edge tags in alternating colors $f$ and $g$ either differing in or adding up to $a$, the latter
occurring precisely if both $f\leq a$ and $g\leq a$.

Denoting a path $H'$ as above by $L(f,g,a)$, we have:
\begin{enumerate}\item
$f=g$ whenever $f=a/2\in\Z$ or $g=a/2\in\Z$; in this case, $d=a/2$ if $d\geq b,c$;
\item
the edges colored $2a$ in $L(a,2a,a)$ are thick.
\end{enumerate}
If two such paths are
parallel and contiguous in $H'$ then they are expressible as $L(f,g,a)$ and $L(h,\!f,a)$, with $|g-h|=2a$ or $g+h=2a$, the latter
occurring precisely if both $g\leq 2a$ and $h\leq 2a$. Here, $g,h$ are the edge colors opposite in the butterflies taking place between $L(f,g,a)$ and $L(h,\!f,a)$. The edges of $L(f,g,a)$ and $L(h,\!f,a)$ colored with $f$ are divided into pairs of opposite edges of the CHs lying between $L(f,g,a)$ and $L(h,\!f,a)$.

\begin{obs}\label{capis}
Given a vertex $v$ of $H(t,a)$, let $f,g,h,i$ be the colors of the
edges incident to an unfolding vertex of $v$ in $H'(t,a)$. If $a$ is
odd or if $v$ is not in an $L(a/2,a/2,a)$ then there is exactly
one other vertex $u$ of $H$ such that the edges incident to any
unfolding vertex of $u$ in $H'$ have colors $f,g,h,i$. In this case
$u$ and $v$ belong to $s=f\!ghi$ and the edge $uv$ has
color $a$.\end{obs}

We may assume that $v$ is shared in $H(t,a)$ by
$a.f\!gj$ and by $a.hij$ so that the edge of $s$ having $v$ as an
endvertex but not having $u$ as an endvertex is colored with $j$, and $j$ colors $v$ in
$s$.

\section{$K_4$-types of $\Z_n$}\label{s9}

\begin{prop}\label{phi}
Let $0<n=2k+1\in\Z$. There is a colored supergraph $G_{n,4}$
of the graph $G'_{n,4}$ introduced in {\rm Section~\ref{s4}}
and a well-defined transformation $\Phi_n$ from $G_{\infty,4}$
onto $G_{n,4}$ that operates by replacing all colors of $\N$
tagging the objects, e.g. vertices, edges, {\rm CT}s and {\rm CH}s of
$G_{\infty,4}$, by their image colors under reduction {\rm MOD} $n$
in the sense that all vertices $($resp. edges$)$ with a common image {\rm MOD} $n$ color disposition can be identified to a corresponding vertex $($resp. edge$)$.
\end{prop}

\begin{proof} Let $A$ be the subset of vertices of the graph $G_{\infty,4}$ introduced in
Section~\ref{s5}
whose colors have exclusively constituents $\le k$ and let $B$ be the set of neighbors
of vertices of $A$ in $G_{\infty,4}$. Let $F$ be the graph induced by $A\cup B$
in $G_{\infty,4}$.
By reducing MOD $n$ all the colors tagging objects of $F$, the
resulting color identifications in $F$ yield $G_{n,4}$. Note that the reduction MOD $n$ for vertices
happens solely for the vertices of $B$. Once these vertices are reduced MOD $n$,
they have the same colors as some vertices of $A$, so they must be identified
correspondingly, and the edges from $A$ to $B$ are then transformed into edges
joining vertices of $A$ which were not originally induced by $A$ in $G_{\infty,4}$.
Now, $\Phi_n$ is defined by replacing the colors of the
objects in $G_{\infty,4}$ (vertices, edges, CTs and CHs)
by their reductions MOD $n$, which yields the corresponding
objects in $G_{n,4}$. \end{proof}

\begin{obs}\label{xxx}
The graph $G_{n,4}$ is an edge-disjoint union of possibly degenerate {\rm CT}s, at most three incident to each vertex.
\end{obs}

\begin{cor}
$G_{n,4}$ is connected, for any odd positive integer $n$.
\end{cor}

\begin{proof} Apply Corollary~\ref{conn} and Proposition~\ref{phi} to the
(continuous) map $\Phi_n:G_{\infty,4}\rightarrow G_{n,4}$. \end{proof}

Application of $\Phi_n$ to the charts of $G_{\infty,4}$ yields charts of
$G_{n,4}$. The collection of charts of $G_{n,4}$,
($G_{\infty,4}$), whose CT centers are colored $i$, for each $i\in
\{1,\ldots,n/2\}$, is called an {\it $i$-atlas}.

\begin{cor}\label{ult}
Let $\rho_n:[k]\rightarrow\{\mbox{atlases of }G_{n,4}\}$ be the assignment
given by $\rho_n(i)=$ $i$-atlas of $G_{n,4}$, for each $i\in[k]$.
If $\gcd(n,i)=1<i<n/2$, then $\rho_n(i)$ is obtained from
$\rho_n(1)$ by replacing each color $c$ tagging a vertex, edge,
{\rm CT} or {\rm CH} of $\rho_n(1)$ by the reduction {\rm MOD} $n$ of $c.i$. If $n$ is
prime, applying $\Phi_n$ to the $i$-atlases of $G_{\infty,4}$
yields $\lfloor n/2 \rfloor$ $i$-atlases of $G_{n,4}$, which are graph isomorphic.
\end{cor}

\begin{proof}
The given reduction MOD $n$ identifies oppositely signed colors mod $n$. \end{proof}

Chart $\rho_{13}(1)$, depicted in Figure 6 (where a superposition of part of the
$\{6,3\}$-regular hexagonal tessellation ${\mathcal H}$ with its edges intersecting at $90\deg$ some of the edges of $\rho_{13}(1)$ is shown in relation to Figure 7 below) exemplify the following properties, which follow by combining the images of the subgraphs 1[011], 2[011], 3[112] under the isomorphisms $\rho_n(1) \rightarrow \rho_n(i)$:

\begin{enumerate}
\item Chart $\rho_n(1)$ is representable in a plane triangle $T(n,1)$
whose sides are SAs of the subgraph $1[011]\subset G_{\infty,4}$, namely two
SAs of type $(2)$ and one of type $(1)$, as in Subsection~\ref{sym}.
\item The internal angle between the SAs of type $(2)$ is $60^\circ$. The internal angles
between each of these and the SA of type $(1)$ are $30^\circ$ and $90^\circ$.
The angle of $30^\circ$ has its vertex at the center $v$ of the CH $1.000$ so $\rho_n(1)$ is represented as a twelfth part of the total angle of $360^\circ$ at $v$.
The angle of $90^\circ$ has its vertex at $0jj1jj,$ where $j=(n-1)/2$.
\item There is only one maximal path $L_{n,1}$ of $\rho_n(1)$ passing
through $0jj1jj$ with its edges having color $j$ and cutting the
opposite side of $T(n,1)$ at $90^\circ$ on a thick edge.
\item The angle of $60^\circ$ has its vertex at the center of the CT $1hhh$,
where $h = (n-5)/2$.
\end{enumerate}

\begin{figure}[htp]
\includegraphics[scale=0.37]{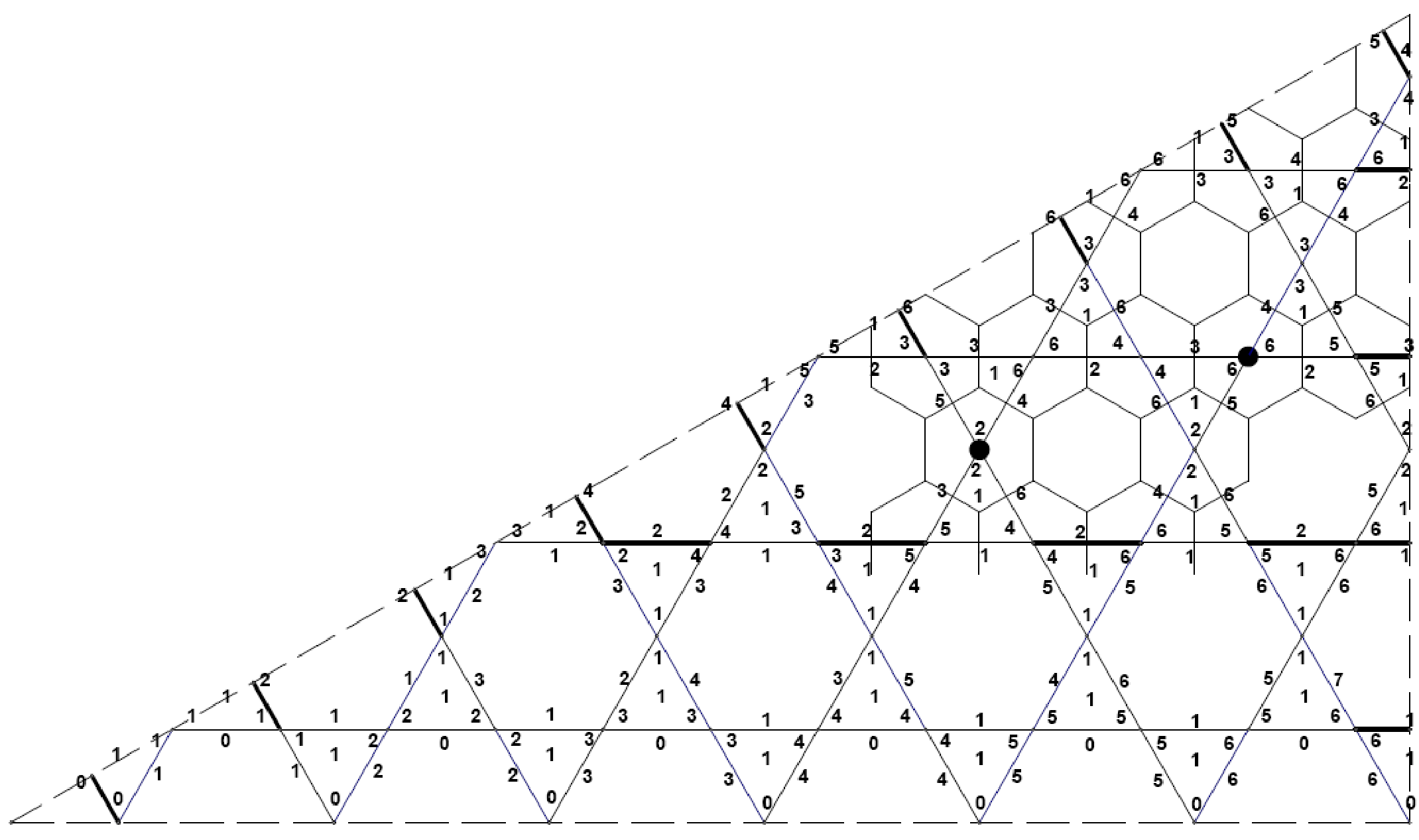}
\caption{Superposition of drawings for $\sigma_n(1)$ and
$\tau_n(1)$}\end{figure}

\begin{prop}\label{23}
The diameter of $G_{n,4}$ is  both $\Omega(n)$ and
$O(|V(G_{n,4})|^{1/3})$, so that the asymptotic diameter of $G_{n,4}$ is
$|V(G_{n,4})|^{1/3}$.\end{prop}

\begin{proof} First, we claim that $|V(G_{n,3})|=O(n\phi(n))$, where
$\phi(n)=$ Euler characteristic of $n$. Every $aa0$, where
$\gcd(a,n)=1$, belongs to $G_{n,3}$. Thus, there are $\lfloor
\phi(n)/2 - 1 \rfloor$ paths whose ends are $011$ and $0aa$, with
$0<a\leq\lfloor n/2 \rfloor$ and $\gcd(a,n) = 1$. But the distance
from $0aa$ to $011$ in $G_{n,3}$ is no more than $a$, yielding our claim.
If we fix a
$K_3$-type of $abcde\!f\in G_{n,4}$, say $abc$, then for each
color $d$ MOD $n$ there are at most two different values for $e$ but
a unique value for $f$. This way, there are at most
$n\phi(n)(2\lfloor n/2\rfloor)$ different $K_4$-types MOD $n$. Thus, $|V(G_{n,4})|=O(n^2\phi(n))$.
Let us see now that the diameter of $G_{n,4}$
is $\Omega(n)$. A path of length $n+1$ between
$110110$ and $112(n-1)nn$ happens along the image of $L(1,2,2)$.
Thus, the diameter of $G_{n,4}$ is both $\Omega(n)$ and
$O(|V(G_{n,4})|^{1/3})$. \end{proof}

A representation of the charts of $G'_{n,4}$ leading to the
connectedness of $G'_{n,4}$ for $n$ large is introduced. Let
$\sigma_n(1)$ be the restriction of $\rho_n(1)$ induced by the rainbow
$K_4$-types. We superpose the T-subgraph representation of
$\sigma_n(1)$ with a $\{6,3\}$-regular hexagonal tessellation
${\mathcal H}=\tau_n(1)$ (\cite{Fej}, page 43) such that: {\bf (a)} each edge
$\epsilon$ of $\sigma_n(1)$ is traversed by an edge $\epsilon'$ of
$\tau_n(1)$ at $90^\circ$ at the common midpoint of $\epsilon$ and
$\epsilon'$; {\bf (b)} each CH of $\sigma_n(1)$ contains in its
interior a regular hexagon of $\tau_n(1)$. Figure 6 contains a
superposition of a representation of $\sigma_{13}(1)$, with the two
rainbow $K_4$-types indicated as bullets $\bullet$ and the part of $\tau_{13}(1)$ used to
represent $\sigma_{13}(1)$ in Figure 7.

\begin{figure}[htp]
\includegraphics[scale=0.368]{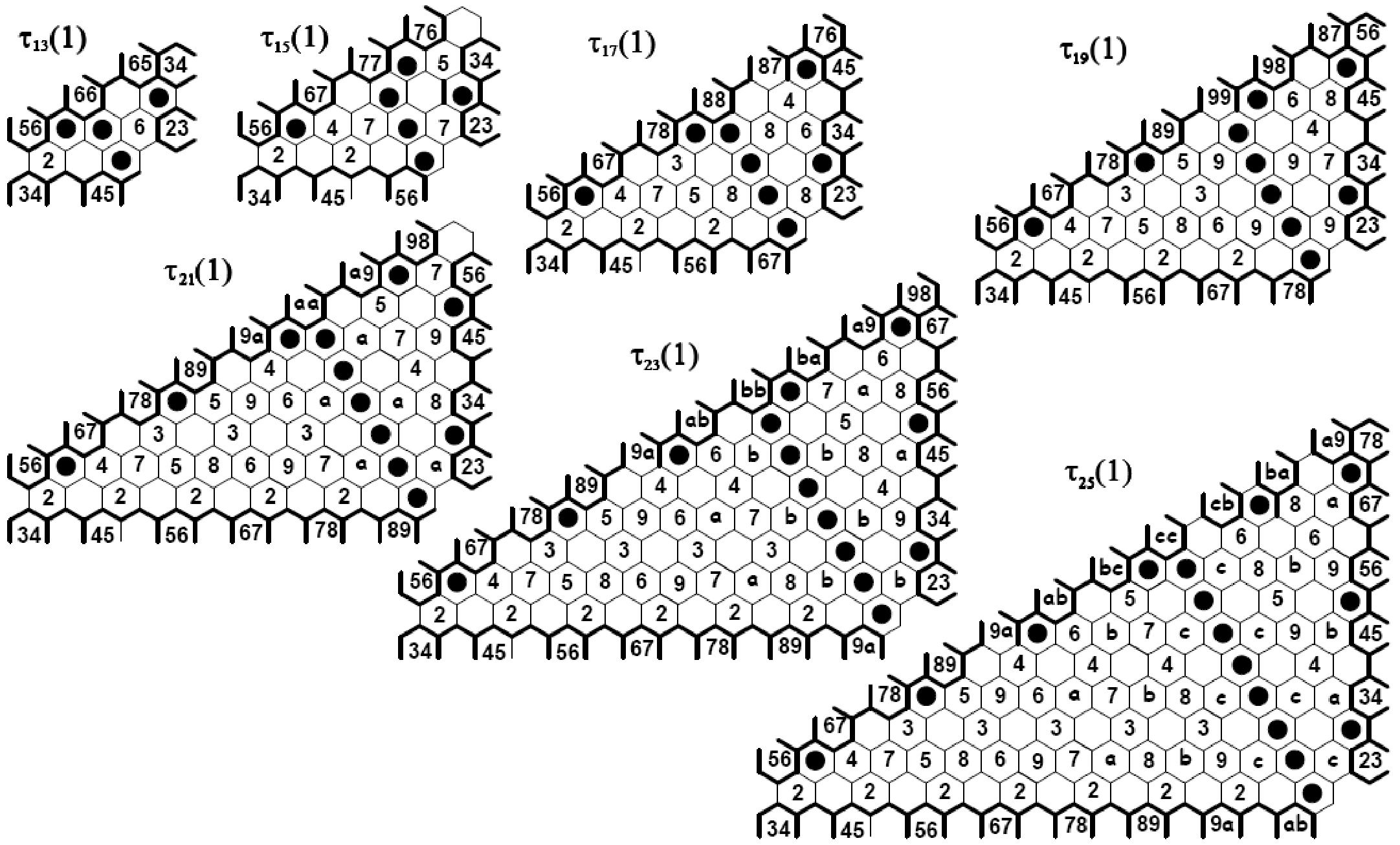}
\caption{The representations $\tau_n(1)$, for $n=13,\ldots,25$}
\end{figure}

In Figure 7, representing $\tau_n(1)$ for odd $n=13,\ldots,25$, each rainbow $K_4$-type of $\sigma_n(1)$ is given  by an hexagon of $\tau_n(1)$ tagged by a positive integer, as suggested in Figure 6 for $n=13$ by the indicated superposition. Each tagged hexagon representing a vertex of $\sigma_n(1)$ is the intersection of two tagged-hexagon sequences in $\tau_n(1)$. There are three directions of parallelism for existing tagged-hexagon sequences: one horizontal and the other two at angles of $\pm 60^\circ$ from the horizontal. Each such sequence is headed on the boundary of $\tau_n(1)$ by a partially-drawn thick-trace hexagon tagged by a pair of integers. Assume the integer tagging an hexagon $\zeta$ of $\tau_n(1)$ is $i$ and the integer pairs heading its two tagged-hexagon sequences are $(p,q)$ and $(r,s)$. Then the $K_3$-types composing $\zeta$ are: $1pq$, $1rs$ and either $ipr$ and $iqs$ or $ips$ and $iqr$. Here, an hexagon is tagged with a bullet $\bullet$ instead of an integer if it represents a non-rainbow $K_4$-type. Each remaining (non-tagged) hexagon stands for a corresponding CH. It follows that each $\sigma_n(1)$ has at least two isolated vertices, represented in $\tau_n(1)$ by: {\bf (1)} the hexagon tagging 2 at the lower-left corner of
$\tau_n(1)$ (that is the $K_4$-type $134265$); {\bf (2)} the hexagon tagged by $\lfloor n/2 \rfloor$, at the lower-right corner of $\tau_n(1)$ (that is the $K_4$-type $123k(k-2)(k-1)$, where $n=2k+1$). If $n\neq 0$ mod 3 then these are the only two isolated vertices of $\sigma_n(1)$. Otherwise, there is exactly one more isolated vertex in $\sigma_n(1)$ and this is determined by the hexagon tagged by $n/3$ at the upper-right corner of $\tau_n(1)$ (that is the $K_4$-type $1(k-2)(k-1)k(k+1)(k+2)$.)

For $n\geq 17$, the isolated vertices of $\sigma_n(1)$ are
nonisolated in the remaining charts $\sigma_n(i)$, where $i\neq 1$
ranges over the units MOD $n$ from 2 to $\lfloor n/2 \rfloor$.
This suggests the following conjecture.

\begin{con}
$G'_{n,4}$ is a connected graph, for $n\geq 17$.\end{con}

The six charts $\tau_{13}(i)$, for
$i=1,\ldots,6$, represent the same pair of isolated vertices shown
in Figure 2(a$_1$) and 2(a$_2$), which are thus the only components of
$G'_{13,4}$. In addition, the four charts $\tau_{15}(i)$, for
$i=1,2,4,7$, represent only a CT and four isolated vertices.

\section{Proofs of the main results}\label{s10}

\begin{proof} {\bf(of Theorem~\ref{1})}
By Proposition~\ref{23}, the asymptotic diameter of $G_{n,4}$ is
$|V(G_{n,4})|^{1/3}$.
The vertices $v\in V_6$ in any member $G=G_{n,4}$ of ${\mathcal G}_1$
are the rainbow $K_4$-types in $G$. The four $K_3$-types of each such
rainbow $K_4$-type form three distinct pairs of $K_3$-types, each corresponding to a respective triangle of $G$. This yields three
triangles $T_0,T_1,T_2$ almost always distinct as in the statement, so that each pair $\{T_i,T_j\}$ with $i\ne j$ determines two different butterflies at $v$ and
respective charts $D_{i,j}^0$ and $D_{i,j}^1$. Let $S\subseteq V_6$ be composed by these vertices $v$. Clearly, $|S|$ is asymptotically $|V_6|$. Now, $V(G)\setminus V_6$ has its vertices at distance no more than $2$ both from the boundary of charts $\tau_n(i)$ and from the diagonal paths $\eta(i)$ in them,
with these paths departing from boundary vertices realizing angles of $90^\circ$ as in the upper right representation in Figure 5 and as in Figure 6. This insures that $|V(G)\setminus V_6|$ grows linearly as $n$ increases, while $|V_6|$ has a quadratic growth with respect to $n$,
so $V_6$ has asymptotic order $|V(G)|$.
Each of the four $K_3$-types composing the $K_4$-type associated with a vertex of $S$ offers three positive integers that color the edges of a corresponding chart modeled on ${\mathcal H}$ as in \cite{De4}, Theorem 2. Each of these three integers colors the edges of a
parallel class of edges in that chart. These completes the proof of Theorem~\ref{1}.\end{proof}

\begin{proof} {\bf(of Corollary~\ref{2})}
Let ${\mathcal G}'_1\subset{\mathcal G}_1$ be formed by the $G_{n,4}$ with $n$ an odd prime.
Then, the charts $\tau_n(i)$ are pairwise isomorphic. They are related with the graphs $D_{i,j}^k$ as follows, for $i=1,\ldots,\frac{n}{2}$. Each $\tau_n(i)$ has
two components formed by vertices representing rainbow $K_4$-types.
These components are: {\bf(a)} contained in a
$30^\circ$-$60^\circ$-$90^\circ$ triangle $R$ (formed by the three
delimiting SAs); {\bf(b)} separated by the path $\eta(i)$ in
$\tau_n(i)$. The union of the two $30^\circ$-$60^\circ$-$90^\circ$
triangles delimited by the SAs and $\eta(i)$ yields $\tau_n(i)$.
By Corollary~\ref{ult}, there are $\lfloor n/2\rfloor$ charts $\tau_n(i)$. We consider stripping bands of the delimiting
SAs in the $30^\circ$-$60^\circ$-$90^\circ$ triangles in order to get rid of loops. This
reduces the resulting $(n-1)$ $30^\circ$-$60^\circ$-$90^\circ$
triangles. The stripped triangles are split into two halves by the paths $\eta(i)$, each half leading to isomorphic $\mathcal D$-modeled subgraphs, with the vertex numbers in the two halves, for $y\ge 1$, equal to:
$|V'^-_y|=2\sum_{i=1}^{y}i$ and $|V'^+_y|=-2+6\sum_{i=1}^{y}i$, if $k=5+2y$;
resp. $|U'^-_y|=|V'^-_y|-y$ and $|U'^-_y|=|V'^+_y|-3y$, if $k=4+2y$.
By removing from $U'^\pm_y$ (resp. $V'^\pm_y$) the isolated vertices in lower-left (resp. lower-; upper-right) corners in the $\tau_n(1)$ in Figure 7, tagged 2 (resp. $k$; $n/3$ if $n\equiv 0$ mod 3), a maximal connected $\mathcal D$-modeled subgraph $U^\pm_y$ (resp. $V^\pm_y$) is obtained.\end{proof}

\end{document}